\documentclass[11pt,a4paper]{amsart}
\usepackage[T1]{fontenc}
\usepackage{lmodern,amsmath,amssymb,microtype}
\usepackage[margin=26mm]{geometry}
\usepackage[hidelinks]{hyperref}
\hypersetup{pdftitle={The maximum number of intersections of two pseudolinear polygons}}

\newcommand{\head}[1]{\par\smallskip\noindent\textit{#1.}\enspace}

\newtheorem{theorem}{Theorem}
\newtheorem{lemma}[theorem]{Lemma}
\newcommand{\intr}{\operatorname{int}}
\title[Intersections of two pseudolinear polygons]{The maximum number of intersections of two pseudolinear polygons}

\newif\ifanonymous
\anonymousfalse
\ifanonymous
  \author{Anonymous}
\else
  \author[Egor Bakaev]{Egor Bakaev}
  \address{Department of Computer Science, University of Copenhagen}
\fi

\begin{document}

\begin{abstract}
In 1993, Dillencourt, Mount, and Saalfeld conjectured that,
for two simple polygons in general position with odd numbers
$m$ and $n$ of sides, the maximum number of intersections
between their boundaries is $(m-1)(n-1)+2$.
In 2026, Seo proved the conjecture using computational
verification, while Ackerman gave a human-readable proof.

We settle a generalization of this conjecture proposed in 2022
by Ackerman, Keszegh, and Rote for pseudolinear polygons.
In this setting, the sides of the polygons are supported
by pseudolines in a common simple arrangement rather
than by straight lines.
\end{abstract}

\maketitle

\section{Introduction}

In 1993, Dillencourt, Mount, and Saalfeld
\cite{dillencourt1993}  
asked:
\emph{What is the maximum number of boundary intersections of two simple
polygons with $m$ and $n$ sides in a general position?}
The case when $m$ or $n$ is even is easy, but the case when both $m$ and $n$ are odd turned out to be difficult.
The proposed upper-bound proof in \cite{dillencourt1993}  contained an error. 
{\v C}ern{\'y}, K{\'a}ra, Kr{\'a}l', Podbrdsk{\'y}, Sot{\'a}kov{\'a},  and {\v S}{\'a}mal~\cite{cerny2003} 
settled the case of a pentagon and proved the upper
bound $mn-m-\lceil n/6\rceil$ for odd $m\ge n\ge7$.
Ackerman, Keszegh, and Rote~\cite{ackerman2022} proved the bound $mn-m-n+C$ for a large absolute constant $C$.

Recently, Seo~\cite{seo2026} obtained the sharp bound using computational verification, and Ackerman~\cite{ackerman2026} subsequently gave a
human-readable proof. 

\begin{theorem}[\cite{seo2026, ackerman2026}]
\label{thm:polygons}
For odd $m,n\ge3$, the maximum number of boundary intersections of two simple
polygons with $m$ and $n$ sides in a general position is $(m-1)(n-1)+2$.
\end{theorem}

A key ingredient in both proofs is the following lemma of
Ackerman, Keszegh, and Rote.
A \emph{corner} consists of two consecutive sides of a polygon,
and its \emph{wedge} is the convex angular region spanned by
these sides at their common vertex.
Two corners, one from each polygon, form a \emph{gadget} if
each side crosses exactly one side of the other corner and
at least one of the two vertices lies inside the wedge of
the other corner.
A gadget has type $P$ if the vertex of the corner of $Q$
lies inside the wedge of the corner of $P$, and type $Q$
if the reverse holds.

\begin{lemma}[Eight-side lemma, {\cite[Lemma~7]{ackerman2022}}]
\label{lem:eight-straight}
Let $P$ and $Q$ be two simple polygons in
general position.
Suppose two gadgets consist of eight distinct sides and
each side of either gadget crosses both sides of the other
polygon in the other gadget.
Then each gadget has exactly one type, and their types
are different.
\end{lemma}

Ackerman, Keszegh, and Rote~\cite[Section~10]{ackerman2022}
also proposed extending their results and conjectures to
pseudolinear polygons.
An \emph{affine pseudoline} is a properly embedded copy of $\mathbb R$ in the
plane. Any two members of a \emph{simple arrangement} meet in exactly one
point and cross there. No three members meet. We consider two simple closed
polygons $P$ and $Q$ whose sides are compact arcs on distinct pseudolines of
one such arrangement. We call such polygons \emph{pseudolinear}. Every
intersection of their boundaries is a proper crossing in the interiors of two
sides, since a polygon vertex already lies on two supporting pseudolines.

We resolve the conjecture \cite[Section~10]{ackerman2022}, which generalizes Theorem \ref{thm:polygons}.

\begin{theorem}\label{thm:main}
For odd $m,n\ge3$, the maximum number of boundary intersections of two simple
pseudolinear polygons with $m$ and $n$ sides in a common simple affine
arrangement is $(m-1)(n-1)+2$.
\end{theorem}

The proof relies on a pseudolinear version of the eight-side lemma.
The authors of \cite[Section~10.2]{ackerman2022}
formulated it and verified it computationally.
We give a short human-readable proof of this lemma in
Section~\ref{sec:eight}. In Section~\ref{sec:count}, we use
the lemma to establish Theorem~\ref{thm:main} by adapting
the combinatorial argument of Ackerman~\cite{ackerman2026}
to the pseudolinear setting.

\section{The eight-side lemma}\label{sec:eight}

Color $P$ blue and $Q$ red. The pseudoline containing a side is its
\emph{support}. A \emph{corner} consists of two consecutive sides, its
\emph{arms}, and their common vertex. Each arm determines a ray
of its support starting at the vertex. For each support take the closed
pseudohalfplane that contains the other arm. Their intersection is the
\emph{wedge} of the corner, bounded by the two rays. Every
other support meets each pseudohalfplane in a ray, so its intersection with
the wedge is an interval, possibly empty or unbounded. In particular, an
arc on such a support with one endpoint inside the wedge cannot cross both
boundary rays. Throughout, \emph{inside} means in the interior. This wedge
property and the opposite-wedge observation below extend
\cite[Observation~2.3]{ackerman2026} to pseudolines.

A blue corner and a red corner form a \emph{gadget} if each of their four
arms crosses exactly one arm of the other color and at least one of the two
vertices lies inside the other corner's wedge. The gadget has \emph{type
$P$} if the red vertex lies inside the blue wedge, and \emph{type $Q$} if
the blue vertex lies inside the red wedge. A gadget of just one type is
\emph{pure}, and a gadget of both types is \emph{double}. Its two disjoint
blue--red pairs are its \emph{bases}. These are Ackerman's gadgets, bases,
and hooking types~\cite[Section~2]{ackerman2026}. 

The wedge obtained by reversing both rays is called the \emph{opposite
wedge}. We will also use the following observation: a pseudoline crossing both
boundary rays of a wedge avoids its opposite wedge, which lies on the
same side of that pseudoline as the apex. Indeed, between the two crossings
the pseudoline is inside the wedge, and its two remaining rays lie in
the two adjacent sectors, by the single-crossing property.

We use Levi's enlargement lemma~\cite[Theorem~5.1.1]{felsner2017} in the
projective completion of the affine arrangement. In the one-point
compactification, the supports form an embedded finite graph. A small
topological vertex-neighborhood of infinity meets each incident end once.
Its complementary disk contains all polygon sides and crossings and meets
each support twice. The ends of any two supports alternate on the boundary,
so the boundary can be parametrized with the two ends of each support
antipodal. Identifying antipodes gives a projective arrangement including
the image of the boundary as its line at infinity. We work in the affine
chart complementary to this line. Levi's lemma adds a pseudoline through
any two points on no common member of the arrangement, which need not be
simple.

\head{Axis observation (cf.~\cite[Observation~6]{ackerman2022} and
\cite[Observation~2.5]{ackerman2026})}
For any gadget, Levi's lemma supplies a pseudoline $d$ through its two
vertices. It meets each of its four arms only at the corresponding vertex.
If the gadget has type $P$, the segment of $d$ from the blue vertex to the red vertex
starts inside the blue wedge and $d$ cannot meet either boundary
support again. Hence $d$ separates the two blue arms. The red
arms meet different blue arms and hence also lie on opposite sides of $d$.
The same holds by symmetry for type $Q$. Thus at each vertex the two rays
of $d$ lie in its wedge and its opposite wedge. In a pure gadget the
vertex outside the other wedge therefore lies inside its opposite wedge.

The next statement is the pseudolinear analogue of
\cite[Lemma~7]{ackerman2022}, in the local form of
\cite[Lemma~2.7]{ackerman2026}. We prove it by a new topological
argument using the two auxiliary axes.

\begin{lemma}[Pseudolinear eight-side lemma]\label{lem:eight}
Suppose two gadgets use eight distinct sides. If each side of either
gadget crosses both sides of the other color in the other gadget, then
the gadgets are pure and have different types.
\end{lemma}

\begin{proof}
By symmetry it suffices to exclude two gadgets of type $P$. For $i=1,2$,
write $u_i,v_i$ for their blue and red vertices and $B_i,R_i$ for their
wedges. Thus $v_i\in\intr B_i$. Shorten the arms just beyond their last
required crossings so that arms of the same color in different gadgets
are disjoint.

\head{Locating the blue vertices}
Add an axis $d_1$ for gadget~1 as above. Each arm of gadget~2 crosses $d_1$
between its crossings with the two arms of the opposite color. In
particular $d_1$ avoids $u_2,v_2$, so Levi's lemma applies again in the
enlarged arrangement to give an axis $d_2$ for gadget~2. Each arm of
gadget~1 similarly crosses $d_2$.

Let $o=d_1\cap d_2$. The two crossings of $d_1$ with the blue arms at
$u_2$ lie on opposite sides of $d_2$, so $o$ is finite and lies between
them on $d_1$, inside $B_2$. The same argument puts $o$ inside all four
wedges. In the sector bounded by the rays $ou_1$ and $ou_2$, a portion of one
blue arm from each $u_i$ joins the two rays (the arms at $u_i$ meet the
other axis on opposite sides of $o$). These disjoint crosscuts have their
endpoints in the same order from $o$ on both rays. Exchange indices
if necessary so that $u_2$ is the nearer endpoint on its ray. Then $u_2$
lies between $o$ and the crossing of $d_2$ with a blue arm of gadget~1,
hence inside $B_1$.

By the wedge observation, $v_2$ is outside $B_1$, since its red arms
cross both blue arms there. Both $o$ and $v_2$ are on the ray of $d_2$
from $u_2$ inside $B_2$. The interval $d_2\cap B_1$ contains $u_2,o$ but not $v_2$,
so their order is $u_2,o,v_2$. The ray from $v_2$ through $o$
lies inside $R_2$, since $d_2$ meets its boundary supports only at $v_2$.
Thus $u_2$ is inside $R_2$. Put $K=B_1\cap R_2$. 

\head{Locating the red vertices}
The red arms at $v_2$ give two disjoint crosscuts of $B_1$ joining its
boundary rays. They have the same order on both rays. Call the nearer
arm $r$ and the farther arm $s$. The quadrilateral between them is $K$,
since each red halfplane defining $R_2$ contains the other red cut.
All four boundary arcs of $K$ are contained in the actual polygon sides.
Let $e,f$ be the blue arms from $u_2$ crossing $r,s$, respectively. They
start inside $K$ and leave through their prescribed red sides, since they
cannot cross the blue boundary sides or the red sides they miss.

After crossing $r$, the arm $e$ is trapped in the near triangle bounded
by $r$ and the blue sides of gadget~1. Before crossing $s$, the arm $f$ stays in $K$; afterwards it
stays on the side of the support of $s$ not containing $u_1$. This side avoids
the opposite wedge of $B_1$. Consequently $f$ cannot cross both blue
supports: after crossing both it would be in that opposite wedge.
Thus both $e$ and $f$ stay in the same open pseudohalfplane of at least
one of the blue supports at $u_1$.

The red arm from $v_1$ crossing this blue side stays inside $B_1$ until
that crossing, and afterwards cannot meet $e$ or $f$. Its initial
portion therefore crosses both. It avoids $r,s$ and lies in one of the
three regions into which they divide $B_1$. The near region misses $f$,
and the far region misses $e$. Hence $v_1$ lies inside $K$.

\head{The trapped endpoint}
Since both red arms from $v_2$ cross both boundary rays of $B_1$, its
vertex is not in the opposite wedge of $B_1$. Thus $v_2$ is in an
adjacent sector, and $r,s$ enter $B_1$ through the same blue side $b$.
The red portions from $v_2$ to $b$ and the intervening portion of $b$
bound a Jordan triangle $T$. The crossings with $b$ precede the respective crossings with $e,f$
along the red arms from $v_2$, since $e,f$ leave $K$ through $r,s$
respectively. The latter crossings lie on the boundary of $B_2$.
Since $v_2$ is inside $B_2$, the interval property puts both crossings
with $b$ inside $B_2$,
and then puts all of $\partial T$ inside $B_2$. The complement of the
open wedge $\intr B_2$ is connected and unbounded, so all of $T$ lies
inside $B_2$.

The red arm from $v_1$ that crosses $b$ starts inside $K$ and cannot
meet $r,s$ or the other blue side of gadget~1. It leaves $K$ across the
edge of $T$ on $b$ and enters $T$. It cannot leave $T$, since it cannot
meet $r,s$ or cross $b$ again. Its free endpoint lies inside $B_2$,
yet it crosses both $e$ and $f$, contradicting the wedge observation.
\end{proof}

\section{The counting argument}\label{sec:count}

Assume that $m,n$ are odd. Let $G$ be the bipartite
\emph{disjointness graph}. Its vertices are the
sides of $P$ and $Q$, and its edges are the disjoint blue--red pairs.
In particular, sides of opposite colors in different components always
cross. Write $c$ for the number of components and $\mu$ for the cycle
rank of $G$, as in~\cite[Section~2]{ackerman2026}. Thus $\mu=|E(G)|-(m+n)+c$, the number of edges outside a
spanning forest. The number of crossings is
\begin{equation}\label{eq:crossings}
 X=mn-|E(G)|=mn-m-n+c-\mu.
\end{equation}
As in~\cite[Proposition~2.2]{ackerman2026}, the Jordan curve $Q$ alternately enters
and leaves the interior of $P$, so $X$ is even. As $m,n$ are odd,
$c-\mu$ is odd. The desired upper bound is equivalent to $c-\mu\le3$.

\subsection{Transitions between components}

The next lemma is the consecutive-side form of the parity argument of
Ackerman, Keszegh, and Rote~\cite[Proposition~4 and Lemma~5]{ackerman2022},
with the hypothesis used by Ackerman~\cite[Lemma~2.4]{ackerman2026}.
We give a short sign-alternation proof valid for pseudolines.

\begin{lemma}[Odd-corner lemma]\label{lem:odd}
If every red side crosses at least one of two consecutive blue sides
$a,b$, then $a,b$ are the blue arms of a gadget. The same holds with
the colors exchanged.
\end{lemma}

\begin{proof}
Regard the supports of $a,b$ as topological coordinate axes, with $a,b$
on the positive horizontal and vertical rays. Write $x(v),y(v)$ for the
coordinate signs at a red vertex $v$. Assign $v$ the sign $y(v)$ if an
incident red side crosses $a$, and $-x(v)$ if it crosses $b$. In quadrants
II and IV the two choices agree. In I and III an incident red side
crosses exactly one positive arm, by the wedge and opposite-wedge
observations. Unless some vertex in I or III has incident sides crossing
different blue arms, the assignment is therefore well-defined.
Each red side reverses the sign, since it crosses the support of
whichever blue arm is used for the sign assignment exactly once.
This contradicts the odd number of red sides.

At such a vertex in I we have a type $P$ gadget. At one in III, the
red arm crossing $b$ crosses the negative horizontal ray before
meeting $b$. That
crossing and the origin are on the same side of the other red support,
which crosses the positive horizontal ray. Interchanging axes gives
the second halfplane condition. Thus the origin is inside the red
wedge, giving type $Q$.
\end{proof}

Every side misses some side of the opposite polygon: its support is
crossed an even number of times by that odd polygon. Thus every component
of $G$ contains both colors. Form the \emph{transition graph} $H$ on the components by joining
$U$ and $V$ whenever consecutive sides of $P$ lie in these two distinct
components. Those sides have no common neighbor in $G$. By
Lemma~\ref{lem:odd} they belong to a gadget with one base in $U$ and one
in $V$. Its two red sides give the same transition along $Q$. Exchanging
colors shows that either cyclic order defines the same graph $H$.
Moreover, $H$ is connected because the cyclic order of $P$ visits every
component. Call a gadget with bases in $U$ and $V$ a \emph{gadget on $UV$}.

For each type, take the edges of $H$ that carry a gadget of that type,
together with their endpoints. These are Ackerman's two type graphs
\cite[Section~3]{ackerman2026}. Two edges of the same type must share
a vertex. Otherwise gadgets on them would use four different components,
and Lemma~\ref{lem:eight} would exclude their common type. A nonempty
family of pairwise incident edges in a simple graph forms a star or a
triangle. Indeed, given edges $UV,UW$, any edge avoiding $U$ must be
$VW$, and an edge meeting all three belongs to that triangle.

Call two edges of $G$ \emph{near} if they share an endpoint or an edge
of $G$ joins their endpoints. If gadgets on $UV$ and $UW$, with
$V\ne W$, have a common type, then their bases in $U$ are near. Otherwise
all eight sides are distinct. All mixed pairs between the two gadgets
cross, either because they lie in different components or because no edge
joins the two bases in $U$. This contradicts Lemma~\ref{lem:eight}.
These are the two consequences of the eight-side lemma used in
\cite[Lemmas 3.1 and 3.2]{ackerman2026}.

We also need the pure-gadget conclusion of Ackerman's
Proposition~4.4~\cite{ackerman2026}. We give a shorter topological proof of the needed part of that
proposition.

\begin{lemma}[Shared sides]\label{lem:shared}
For distinct components $U,V,W$, pure gadgets on $UV$ and $UW$ have
vertex-disjoint bases in $U$.
\end{lemma}

\begin{proof}
Suppose the bases share a blue side $p=AB$, with other blue arms
$a,b$ in $V,W$ at $A,B$, respectively. Both red arms of the gadget at
$A$ cross $b$, and both red arms of the gadget at $B$ cross $a$, since
these mixed pairs belong to distinct components. Write $v_1,v_2$ for
the red vertices and $t_1,t_2$ for the red arms crossing $p$.

If the first gadget is pure of type $P$, the axis observation puts
$A$ inside the opposite red wedge. Its crossing with $t_1$ lies on
one bounding ray of the red wedge, so along $p$ the other red support
is crossed before this crossing. The continuation of $p$, including
$B$, is therefore inside the red wedge. But $b$ crosses both red arms,
contrary to the wedge property. The same holds at $B$.
Both gadgets must therefore have type $Q$.

By oddness, some red side misses the entire support $\ell$ of $p$.
It belongs to $U$, so both $a$ and $b$ cross it and lie on the same
side of $\ell$ (apart from their endpoints on $p$). Let $H_a,H_b$ be
the closed pseudohalfplanes of the supports of $a,b$ containing $B,A$,
respectively, and put $W_0=H_a\cap H_b$. The segment $p$ joins the two
boundary rays of this wedge. By the axis observation applied to the
pure gadgets, $v_1$ is outside $H_a$, $v_2$ outside $H_b$, and both
vertices lie on the side of $\ell$ opposite $a,b$.

Along $t_1$ from $v_1$, its crossing with the support of $a$ precedes
its crossing with $p$, which in turn precedes its crossing with $b$.
Thus $t_1$ enters $W_0$ across the boundary on the support of $a$ and
leaves across $b$; symmetrically $t_2$ enters across the support of $b$
and leaves across $a$. On each bounding ray of $W_0$ the two crosscut
endpoints are on opposite sides of $\ell$. Since $\ell$ crosses both
rays, the portions of the rays on one side of $\ell$ are either both
nearer the apex or both farther from it. The endpoints of the two
crosscuts therefore occur in opposite orders from the apex, forcing
their interiors to intersect by Jordan separation. This contradicts
simplicity of $Q$. The case of a shared red side is symmetric.
\end{proof}

\subsection{Cuts and forests}

We use the following reformulation of Ackerman's cut argument
\cite[Proposition~4.5]{ackerman2026}.

\head{Separator rule}
Partition $V(G)-\{z\}$ into two parts, each containing a side of the
color of $z$, with no edge of $G-z$ between the parts. In the linear
polygon order after deleting $z$, two consecutive sides of that color
lie in different parts. They have no common neighbor in $G$, so by
Lemma~\ref{lem:odd} there is a gadget with one base in each part, avoiding $z$.

\head{Cut rule}
For a cut $D\mid E$ of $H$ with $|D|\ge2$ and $E\ne\varnothing$, there
are two gadgets across the cut with disjoint $D$-bases (bases in components
of $D$). Otherwise, since pairwise incident edges of a bipartite graph
have a common endpoint, all $D$-bases across the cut contain a vertex $z$.
Delete $z$ and apply the separator rule to the two parts given by $D,E$:
each still contains a side of the color of $z$, since every component
contains both colors. This is a contradiction.

\head{Forest rule (cf.~\cite[Proposition~4.3]{ackerman2026})}
In a forest, three disjoint pairwise near edges are impossible: contracting
them would give a triangle minor. Also, all edges near both of two
non-near edges contain a common vertex. Indeed, if the two edges lie in
one tree, choose a vertex strictly between them on their unique connecting
path; any edge near both must contain that vertex. If they lie in different
trees, no edge is near both.

The following is a short cycle-rank argument for the capacity bound of
Ackerman~\cite[Proposition~4.7]{ackerman2026}.

\begin{lemma}[Capacity bound]\label{lem:capacity}
At a component $U$, select at most two gadgets for each neighbor in $H$,
with disjoint $U$-bases if two are selected for the same neighbor.
Let $k$ be their number and $\mu(U)$ the cycle rank of $U$. Then
\[
 k\le 2\mu(U)+7,
 \qquad k\le \mu(U)+3\quad\text{if all are pure of one type}.
\]
\end{lemma}

\begin{proof}
Each side occurs in at most two selected bases, since it has only two
cyclic neighbors and occurs at most once per neighbor. If all gadgets
are pure, the bases are pairwise disjoint by Lemma~\ref{lem:shared} and
the selection rule.
Delete one endpoint of each edge outside a spanning tree of $U$, discarding
at most $2\mu(U)$ gadgets (or $\mu(U)$ if all are pure). The remaining
bases are in a forest, and nearness between them is unchanged.

Consider the surviving bases of one type, with incidence at most $d$ at
a vertex, where $d\in\{1,2\}$ ($d=1$ if all surviving gadgets are pure).
If they are pairwise near, they number at
most $2d$: they split into $d$ matchings, each of size at most two by the
forest rule (for $d=2$, alternate the edges of paths and even cycles).
Otherwise two non-near bases must be toward the same neighbor. All others
are near both and hence contain a common vertex, giving at most $2+d$
bases in total. Thus each type contributes at most four, or at most three
if all surviving gadgets are pure. If a double gadget survives it is counted in
both types; otherwise each type contributes at most three. Hence at most
seven gadgets survive, or at most three if all are pure of one type.
\end{proof}

\subsection{The upper bound}

We finish by adapting and simplifying Ackerman's star/triangle analysis
\cite[Section~5]{ackerman2026}.

\begin{proof}[Proof of the upper bound in Theorem~\ref{thm:main}]
We show that $c\le\mu+4$, so assume $c\ge5$. Recall that each type graph
is empty, a star, or a triangle, and their union is connected.

\head{A common star}
If $H$ is a star with center $U$, the cut rule supplies two gadgets per
leaf with disjoint $U$-bases. By Lemma~\ref{lem:capacity},
$2(c-1)\le2\mu(U)+7$, so $\mu\ge\mu(U)\ge c-4$.

\head{At least one star type, but $H$ is not a star}
Both type graphs are nonempty. If one is a triangle, contract its three
nodes into a single group; otherwise take every node as a singleton
group. Let $s$ be the number of star centers ($s=1$ or $2$). There are
$c+s-4$ groups not containing a star center, each adjacent to one or two
actual star centers. For a group adjacent to just one center, the cut
rule supplies two gadgets at that center with disjoint bases there; its
complement has at least two nodes. For a group adjacent to both centers,
choose one gadget at each center. All the selected gadgets are pure of
the respective star types, and selections at each center satisfy the
hypotheses of Lemma~\ref{lem:capacity}: distinct groups involve distinct
neighbors. Summing the pure bound over the centers gives
\[
 2(c+s-4)\le \mu+3s,
 \qquad\text{so}\qquad \mu\ge 2c-8-s\ge c-4,
\]
except when $c=5$ and $s=2$.

In that case suppose $\mu=0$. Each of the two centers $S,T$ has three
selected gadgets by Lemma~\ref{lem:capacity}, since six were selected in all. Neither has three noncentral
neighbors: three pure bases from distinct neighbors would be disjoint
and pairwise near in a tree. Thus each center has one private leaf and
they share the third noncentral node. At $S$, the two bases $h,h'$ toward
its private leaf $L$ are non-near, since the base toward the common
neighbor is disjoint from and near both. The cut rule for $\{S,L\}$
supplies two gadgets with disjoint bases at $S$. They cannot both have
the type of $S$: both bases would be near $h,h'$ and hence contain a
common vertex. Thus one is pure of the other type, necessarily on $ST$.
Adding it at $T$ gives four selected pure gadgets of the same type in a
tree, contradicting Lemma~\ref{lem:capacity}.

\head{Two triangles}
Connectedness forces $c=5$, with exactly one common node $U$. Suppose
$\mu=0$. All transition gadgets are pure, since each transition has a
vertex-disjoint transition in the other triangle. For each triangle,
the endpoints of all its $U$-bases induce a connected subtree of $U$ by
nearness across its two neighbors. The subtrees are vertex-disjoint by
Lemma~\ref{lem:shared}. Comparison with the edge opposite $U$ in each
triangle excludes gadgets with both bases in $U$ by Lemma~\ref{lem:eight}.
Choose a vertex $z$ on the path joining the subtrees, possibly an endpoint,
so that no component of $U-z$ meets both subtrees. Partition $V(G)-\{z\}$
according to the two triangles, assigning components of $U-z$ meeting a
subtree to its triangle's part and the others arbitrarily. Neither an
edge of $G-z$ nor a gadget avoiding $z$ crosses this partition. Both
parts contain sides of the color of $z$, contrary to the separator rule.
Thus $\mu\ge1=c-4$.

In all cases $c-\mu\le4$. Since $c-\mu$ is odd, $c-\mu\le3$;
substitution in~\eqref{eq:crossings} finishes the proof.
\end{proof}

\subsection{The cases with an even number of sides}
If $m$ is even, the maximum is $mn$ for even $n$ and $m(n-1)$ for odd $n$.
Indeed, each pair of sides crosses at most once, and $Q$ crosses the
support of each side of $P$ an even number of times. For odd $n$ at most
$n-1$ sides of $Q$ can therefore cross a given side of $P$.

\subsection{Attainment of the bound}
The straight-sided constructions (e.g.,~\cite{ackerman2022})
attain the stated bounds for all parity combinations. A 
small perturbation make the position general. Hence the constructions
are pseudolinear and the bound in Theorem~\ref{thm:main} is attained.

\subsection*{Aknowledgments.}
We thank Eyal Ackerman for useful remarks.

\subsection*{Funding}
The author is supported by the DNRF Chair grant 20231101-28697.

\subsection*{AI usage.} ChatGPT was used in developing the
proof and preparing the text.

\bibliographystyle{amsalpha}
\bibliography{pseudolinear}

\end{document}